\documentclass[11pt]{article}

\usepackage{amsfonts}
\usepackage{amsthm}
\usepackage{amsmath}
\usepackage{epsfig}
\usepackage{color}
\usepackage{hyperref}
\usepackage[all,tips]{xy}

\DeclareMathOperator{\PSL}{PSL}

\def\qed{ $\sqcup\!\!\!\!\sqcap$}

\newcommand{\set}[1]{\left\{#1\right\}}

\newtheorem{theorem}{Theorem}[section]
\newtheorem{lemma}[theorem]{Lemma}

\newtheorem{conjecture}{Conjecture}

\title{\textbf{Determining hyperbolic \\ 3--manifolds by their surfaces} }
\author{D.~B. McReynolds\thanks{This work was supported in part by the NSF}~~and
~A.~W. Reid\thanks{This work was
supported in part by the NSF}}

\begin{document}
\maketitle

\begin{abstract}
\noindent In this article, we prove that the commensurability class of a closed, orientable, hyperbolic 3--manifold is determined by the surface subgroups of its fundamental group. Moreover, we prove that there can be only finitely many closed, orientable, hyperbolic 3--manifolds that have the same set of surfaces.
\end{abstract}

\section{Introduction}

The geodesic length spectrum of a Riemannian manifold $M$ is a basic invariant that has been well-studied due to its connection with the geometric and analytic structure of $M$. When $M$ has negative sectional curvatures, there is a strong relationship between this spectrum and the eigenvalue spectrum of the Laplace--Beltrami operator (see \cite{DG},\cite{Gan}), and the latter is well known to determine basic geometric/topological invariants like the dimension, volume, and total scalar curvature of $M$. 

In this article, we focus on variations of the surface analog of the geodesic length spectrum for closed, orientable, hyperbolic 3--manifolds introduced by the authors in \cite{MR} (see also \cite{LM}, \cite{Meyer}). We take this theme further and study {\em the full surface spectrum (or set)} of such manifolds (see \S 2 for definitions) which loosely takes into account all of the $\pi_1$--injective surface subgroups of the fundamental group of $M$. Our main result can be informally stated as follows (see Theorem \ref{maincor} for the precise statement).  

Recall that two Riemannian manifolds $M_1,M_2$ are \emph{commensurable} if there exists a Riemannian manifold $M$ and finite Riemannian covers $M \to M_1,M_2$. 
\begin{theorem}\label{main}
For any closed, orientable, hyperbolic $3$--manifold $M$, there are at most finitely many non-isometric closed, orientable, hyperbolic 3--manifolds with the same surface set as $M$. Furthermore, all such manifolds are commensurable.
\end{theorem}

By way of comparison, for the eigenvalue or geodesic length spectra, many commensurability
and finiteness results have been established. The second author
\cite{Re} proved that isospectral (i.e.~the same eigenvalue spectra)
or length isospectral (i.e.~the same geodesic length spectra)
arithmetic hyperbolic 2--manifolds are
commensurable (see \cite{MacR} for a thorough treatment of arithmetic hyperbolic 2-- and 3--manifolds). Chinburg--Hamilton--Long--Reid \cite[Thm 1.1]{CHLR}
proved an identical result for arithmetic hyperbolic
3--manifolds. Prasad--Rapinchuk \cite[Thm 8.12]{PR} determined when
these commensurability rigidity results hold for general,
arithmetic, locally symmetric orbifolds, proving that in many settings
the commensurability class of the manifold is determined by the
eigenvalue or geodesic length spectra. It was already known that the
commensurability class is not always determined by these spectra as
Lubotzky--Samuel--Vishne \cite[Thm 1]{LSV} produced higher rank,
arithmetic, locally symmetric incommensurable isospectral examples
prior to \cite{PR}. In \cite[Thm 1.1]{MR}, the
authors proved a result similar to Theorem \ref{main}. Namely, if
$M_1,M_2$ are arithmetic hyperbolic 3--manifolds that contain a totally geodesic
surface, and have the same set of totally
geodesic surfaces, then they are 
commensurable. Meyer \cite[Thm C]{Meyer} established
a higher dimensional analog for certain classes of arithmetic
hyperbolic $n$--manifolds. It is worth emphasizing that our present
work differs from all the above works in one important and fundamental
way. Namely, we do not impose an arithmetic assumption.

In an effort to see whether or not \cite{Re} holds in the non-arithmetic setting, Millichap \cite{Mill2} constructed $(2n)!$ incommensurable, non-arithmetic hyperbolic 3--manifolds with the same first $2n+1$ (complex) geodesic lengths. The manifolds have the same volume and the volume of these manifolds is linear in $n$. Since the completion of this paper, Futer--Millichap \cite{FuterMill} and Linowitz--McReynolds--Pollack--Thompson \cite{LMPT} have produced additional examples of non-arithmetic and arithmetic hyperbolic 2-- and 3--manifolds that share the same geodesic lengths for the first $n$ lengths or any finite subset of lengths, respectively. Both constructions give control on the volumes of the examples as well. 

In \cite[Thm 1.2]{MR}, examples of non-isometric, closed, hyperbolic 3--manifolds with the same spectra of totally geodesic surfaces were constructed (see also \cite[\S 5]{Mc} and \cite{MMS}). Those methods can be employed to also produce arbitrarily large finite sets of non-isometric closed hyperbolic 3--manifolds $\set{M_j}$ that pairwise have the same totally geodesic surface spectra (the spectra can be ensured to be infinite as well). However, it is unknown if an infinite set of such manifolds can exist. In particular, the totally geodesic surface analog of our finiteness result is unknown. Finally, for the full surface spectrum, there are no known examples of non-isometric hyperbolic 3--manifolds $M_1,M_2$ with the same full surface spectrum. 

\smallskip\smallskip

\noindent \textbf{Question.}
\textsl{Do there exist non-isometric hyperbolic 3--manifolds with the same full surface spectrum?}

\smallskip\smallskip

\noindent{\bf Acknowledgements:}~{\em The authors are very grateful to the referee for their careful reading of the paper, and several useful comments and suggestions.}

\section{Notation and Preliminaries}

Throughout, $M={\bf H}^3/\Gamma$ will be a closed, orientable,
hyperbolic 3--manifold and $\Sigma_g$ will denote the closed
orientable surface of genus $g$. It was proved by Thurston \cite[Cor
8.8.6]{T} that the number of $\Gamma$--conjugacy classes of subgroups
of $\Gamma$ isomorphic to $\pi_1(\Sigma_g)$ is finite. A breakthrough
was provided by Kahn--Markovic \cite[Thm 1.1]{KM} who proved that this number is non-zero for
certain values of $g$. Building on previous work of Masters \cite[Thm 1.2]{M}, Kahn--Markovic \cite[Thm 1.1]{KM1} also provide estimates for these numbers.

For each discrete, faithful representation $\rho\colon \pi_1(\Sigma_g) \to \PSL(2,\bf
C)$, we refer to the image $\Delta_\rho$ as a \emph{Kleinian surface
  group}. For each $\Delta_\rho$, let $\ell_\rho(M)$ denote the number
of $\Gamma$--conjugacy classes of subgroups $\Delta < \Gamma$ that are
$\PSL(2,\bf C)$--conjugate to $\Delta_\rho$. Typically the value of
$\ell_\rho(M)$ will be zero (e.g.~for those $\Delta_\rho$ that contain
an element with transcendental trace). We define the \emph{full surface spectrum} of $M$ 
to be the set 

$${\cal S}(M) = \{(\Delta_\rho,\ell_\rho(M))~:\ell_\rho(M) \ne 0\}.$$
 We define the {\em surface set} of $M$ to be the set $\mathrm{S}(M) =
\{\Delta_\rho~:\ell_\rho(M) \ne 0\}$. 
The case when $\Delta_\rho$ is
Fuchsian was studied in \cite{MR}, and we denote the associated spectrum here by $\mathrm{S}_{Fuc}(M)$ and call this the
\emph{genus spectrum}. In this note, particular emphasis will be placed upon those Kleinian surface groups $\Delta_\rho$ corresponding to virtual fiber subgroups of $\Gamma$. These subgroups arise as fibers of mapping tori of finite covers of $M$. That is, if $\pi_1(S) = \Delta_\rho$, then there is a finite cover $M' \to M$ and a diffeomorphism $\psi\colon S \to S$ such that $M'$ is the mapping torus of $S$ with respect to $\psi$ where $\Delta_\rho$ is the kernel of the associated homomorphism $\pi_1(M') \to \mathbf{Z}$. 

It is a well-known consequence of  the solution to the Tameness Conjecture (by Agol \cite{Ag} and Calegari--Gabai \cite{CG}) and Canary's Covering Theorem \cite{Can}
that these virtual fiber subgroups of $\Gamma$ are precisely those $\Delta_\rho < \Gamma$ that are finitely generated and geometrically infinite (see also the earlier work of Bonahon \cite{Bon} and Thurston \cite{T}).  Since being geometrically infinite depends only on $\rho$, these surface subgroups provide an important subclass of surface subgroups that can be used to control the topology of 3--manifolds. For future reference, we denote the associated spectrum for this class of surface subgroups by
\[ \mathrm{S}_{vf}(M) = \{\Delta_\rho\in \mathrm{S}(M)~:\Delta_\rho~\hbox{is a virtual fiber subgroup}\}. \]
Essential in our work is the groundbreaking work of Agol \cite{Ag1} and the aforementioned work of Kahn--Markovic \cite{KM}. 

We summarize from their collective work what is needed for us in the following theorem.

\begin{theorem}\label{summary}
Let $M={\bf H}^3/\Gamma$ be a closed, orientable, hyperbolic 3--manifold. Then
\begin{itemize}
\item[(a)] $\mathrm{S}(M)\neq \emptyset$.
\item[(b)] $\mathrm{S}_{vf}(M)\neq \emptyset$.
\item[(c)] $\mathrm{S}_{vf}(M)$ contains infinitely many elements $F_\rho$ that
are not commensurable and in particular have arbitrarily large genus.
\end{itemize}
\end{theorem}

\noindent {\bf Proof:}~Given the preamble to the statement of the theorem, only (c) requires comment. By \cite{Ag1} there is a finite sheeted cover $M_0 \rightarrow M$ such that $b_1(M_0)\geq 2$ and $M_0$ is fibered.  In particular, by \cite{T2}, $M_0$ is fibered in infinitely many different ways. Indeed, it follows from \cite{T2} that we can find fibered surfaces of arbitrarily large genus occurring as integral lattice points in the (open) cone over a top dimensional face of the Thurston norm ball.  

Moreover, by \cite[Corollary 3.7]{CSW} we can assume that infinitely many of these provide incommensurable virtual fibers.\qed

\section{Proof of Theorem \ref{main}}

We now state the precise version of Theorem \ref{main} that we will prove in this section.

\begin{theorem}\label{maincor}
If $M$ is a closed, orientable, hyperbolic 3--manifold, then the set
\[ \mathcal{S}_M = \set{N~:\mathcal{S}(M) = \mathcal{S}(N)} \]
is finite. Moreover, if $N \in \mathcal{S}_M$, then $M,N$ are commensurable.
\end{theorem}

As noted above, since being a virtual fiber depends only on $\Delta_\rho$ and not on the ambient manifolds, if $\mathcal{S}(M) = \mathcal{S}(N)$, then $\mathrm{S}_{vf}(M) = \mathrm{S}_{vf}(N)$. In particular, to prove Theorem \ref{maincor}, it suffices to prove the following result.

\begin{theorem}\label{mainprecise}
If $M$ is a closed, orientable, hyperbolic 3--manifold, then the set
\[ \mathcal{S}_{M,vf} = \set{N~:~\mathrm{S}_{vf}(M) = \mathrm{S}_{vf}(N)} \]
is finite. Moreover, if $N \in \mathcal{S}_{M,vf}$, then $M,N$ are commensurable.
\end{theorem}

\noindent {\bf Proof of Theorem \ref{mainprecise}:}~We first prove
that if $\mathrm{S}_{vf}(M) = \mathrm{S}_{vf}(N)$, then $M,N$ are
commensurable. To that end, let $\Delta=\Delta_\rho$ denote a common
virtual fiber subgroup and set $g$ to be the genus of $\Delta$. Since
$\Delta$ is a virtual fiber, we can find pseudo-Anosov maps
$\phi,\psi\colon \Sigma_g \longrightarrow \Sigma_g$ so that
$M_{\phi}\rightarrow M$, $M_{\psi}\rightarrow N$ are finite sheeted
covers and $\pi_1(M_{\phi}), \pi_1(M_{\psi})$ have a common fiber
group $\Delta$. Associated to the fiber group $\Delta$ is a unique
pair of ending laminations $\nu^{\pm}$ in the projective measured
lamination space of $\Sigma_g$ which are left invariant by $\phi,
\psi$ (see \cite{Bon}).  As a result, there exist integers $r,s$ such
that the mapping classes $\phi,\psi$ satisfy $\phi^r =
\psi^s$. Consequently, the bundles $M_{\phi^r}$ and $M_{\psi^s}$ are
isometric. In particular, we have
\[ \xymatrix{ & M_{\phi^r} \cong M_{\psi^s} \ar[rdd]^{\mathrm{finite}} \ar[ldd]_{\mathrm{finite}} & \\ & & \\ M_\phi \ar[dd]_{\mathrm{finite}} & & M_\psi \ar[dd]^{\mathrm{finite}} \\ & & \\ M & & N} \]
and thus conclude that $M,N$ are commensurable. 

It remains to establish the finiteness of $\mathcal{S}_{M,vf}$. 
We will argue by contradiction, and to that end, we assume that there are
infinitely many non-isometric $M_i={\bf H}^3/\Gamma_i$, $i=1,2,
\ldots$ with $\mathrm{S}_{vf}(M)=\mathrm{S}_{vf}(M_i)$ for all $i$. We
will prove that there is $i \geq i_0$ such that the groups $\Gamma_i$ have
uniformly bounded rank. We will then show that for an even larger $i_1$,
the groups $\Gamma_i$ for $i\geq i_1$ must have rank larger than
this uniform bound. Towards that goal, we first assert that the volumes
of the manifolds $M_i$ must be unbounded. Specifically, we have the
following general lemma.

\begin{lemma}\label{VolumeLemma}
The set of volumes for any infinite set $\set{M_i}$ of commensurable, finite volume, hyperbolic $3$--manifolds is unbounded.
\end{lemma}

\noindent \textbf{Proof:} We split into two cases depending on whether
the manifolds are arithmetic or not. Note that since arithmeticity is
a commensurability invariant, either all of the $M_i$ are arithmetic
or all of the $M_i$ are non-arithmetic. If the $M_i$ are arithmetic,
the assertion follows from work of Borel \cite{Bo} since there are
only finitely many arithmetic hyperbolic 3--manifolds of bounded
volume. If the $M_i$ are non-arithmetic, by work of Margulis
\cite{Mar}, there is a unique maximal lattice in the common
commensurability class that contains all of the $\Gamma_i$ as finite
index subgroups. In particular, all the $M_i$ cover the fixed closed
hyperbolic 3--orbifold $Q$ associated to this unique maximal
lattice. Since $Q$ has only finitely many degree $d$ covers for any
$d$, the covering degrees must go to infinity. Consequently, the
volumes cannot be bounded in this case
either. \qed\smallskip\smallskip

Recall that the \emph{injectivity radius at a point $p \in M$} is the 
largest radius for which the exponential map at $p$ is a diffeomorphism. The \emph{injectivity radius of $M$} (which we denote by $\mathrm{InjRad}(M)$) is the infimum of the injectivity radius of $M$ at $p$ over all $p \in M$.

We note that in our setting, as the manifolds $M_i$ are all commensurable, there is also a uniform lower bound of the injectivity radii of the $M_i$.  This again can be established using the arithmetic/non-arithmetic dichotomy. Briefly, in the non-arithmetic case, the injectivity radius is bounded below by the injectivity radius of the orbifold associated to the maximal lattice in the commensurability class containing the $M_i$. In the arithmetic setting, since the injectivity radius is one half of the systole (which is the length of the shortest, closed and non-trivial geodesic), the systole of any orbifold in the commensurability class is uniformly bounded below by a constant that depends only on the invariant trace field (see \cite[Ch 12, p.~378--380]{MacR}). For future reference, we denote the lower bound for $\mathrm{InjRad}(M)$ by $s$. 

Thus we can now assume that we have a sequence of manifolds $M_i$ with injectivity radius at least $s$ and whose volumes get arbitrarily large. We now show how to use this to bound the ranks of the groups $\Gamma_i$ for $i$ sufficiently large. Towards that goal, set $\Delta_0$ to be a common, minimal genus, virtual fiber group in $\Gamma_i$, and set $g$ to be this common, minimal genus.  In order to control the ranks of the groups $\Gamma_i$, we will utilize a quantitative virtual fibering result of Soma \cite[Thm 0.5]{So}.  To state his result, let $\mathrm{Vol}(M)$ denote the volume of $M$, and set $d_1(g,s) = {4s(g-1)\over{\sinh^2(s/2)}}$.

\begin{theorem}[Soma]\label{soma}
If $M$ is a closed, orientable, hyperbolic 3--manifold with 
\[ \mathrm{InjRad}(M) \geq s~\textrm{and}~ \mathrm{Vol}(M) \geq 2\pi d_1(g,s)\sinh^2(d_1(g,s)+1), \] 
then any immersed virtual fiber in $M$ of genus $g$ is embedded.
\end{theorem}

Theorem \ref{soma} with the above conditions on
$\mathrm{InjRad}(M_i),\mathrm{Vol}(M)$ implies that there is $i_{g,s}
\in\mathbf{N}$ such that if $i \geq i_{g,s}$, the virtual fiber group
$\Delta_0$ corresponds to an embedded incompressible surface 
of genus $g$ in $M_i$. This incompressible surface limits the possibilities for the $M_i$. Specifically, $M_i$ is either a
fiber bundle over the circle with fiber group $\Delta_0$, or $M_i$ is
the union of two twisted $I$--bundles. Moreover, in the latter case,
we have a double cover $N_i\rightarrow M_i$ such that $N_i$ is a fiber
bundle with fiber group $\Delta_0$ (see \cite[Ch 11]{He}).

We now leverage the above fiber bundle structure to obtain bounds for the rank of $\Gamma_i$ for $i$ sufficiently large. The rank of $\Gamma_i$ will be denoted by $\mathrm{Rank}(\Gamma_i)$.
\smallskip\smallskip

\begin{lemma}\label{RankLemma}
There exists $i_0 \geq i_{g,s}$ such that if $i \geq i_0$, then $g +1 \leq \mathrm{Rank}(\Gamma_i) \leq 2g+2$.
\end{lemma}

\noindent{\bf Proof:}~We assume throughout that $i \geq i_{g,s}$. Let
$\mathcal{I}_1$ to be the set of $i \geq i_{g,s}$ such that $M_i$ is a
fiber bundle with fiber group $\Delta_0$ and let $\mathcal{I}_2$ to be
the set of $i \geq i_{g,s}$ such that $M_i$ is double covered by $N_i$
where $N_i$ is a fiber bundle with fiber group $\Delta_0$. We first
consider $\set{M_i}_{i \in \mathcal{I}_1}$.  We know from the proof of
the commensurability invariance of $\mathrm{S}_{vf}$ that each $M_i$
must have the form $M_{\phi^{r_i}}$ for some pseudo-Anosov element
$\phi$. Applying Souto \cite[Thm 1]{sou}, there exist $i' \in
\mathbf{N}$ such that $\mathrm{Rank}(\Gamma_i)=2g+1$ for all $i \geq
i'$. Next, we consider $\set{M_i}_{i \in \mathcal{I}_2}$ and apply the
above argument to $N_i$. We obtain $i'' \in \mathbf{N}$ such that
$\mathrm{Rank}(\pi_1(N_i)) = 2g+1$ for all $i \geq i''$. As $N_i$ is a
double cover of $M_i$, we can adjoin one element to $\pi_1(N_i)$ to
generate $\pi_1(M_i)$. Therefore, $\mathrm{Rank}(\Gamma_i) \leq 2g+2$
for all $i \geq i''$. 

Now, set $i_0 = \max\set{i',i''}$ and note that
$\mathrm{Rank}(\Gamma_i) \leq 2g+2$ for all $i \geq i_0$. For the
lower bound, by the Nielsen--Schreier formula, we have
$\mathrm{Rank}(\pi_1(N_i)) \leq 2\mathrm{Rank}(\Gamma_i) - 1$ for all
$i \geq i_0$ and $i \in \mathcal{I}_2$. In particular, $g +1 \leq
\mathrm{Rank}(\Gamma_i)$ for all $i \geq
i_0$. \qed\smallskip\smallskip

We now use Lemma \ref{RankLemma} to complete the proof of Theorem
\ref{mainprecise}. By Theorem \ref{summary} (c), we can find
incommensurable virtual fiber subgroups of arbitrarily large
genus. Choosing a virtual fiber subgroup $\Delta_1$ of genus $g_1$
with $2g+2 < g_1 + 1$ and repeating the above argument, we obtain an
integer $i_1 \geq i_{g_1,s}$ such that $g_1 + 1 \leq
\mathrm{Rank}(\Gamma_i)$ for all $i \geq i_1$. For all $i \geq
\max\set{i_0,i_1}$, we must have $g_1 + 1 \leq \mathrm{Rank}(\Gamma_i)
\leq 2g+ 2 < g_1 +1$, a contradiction. Hence $\mathcal{S}_{M,vf}$ is
finite as required.\qed\\[\baselineskip]

\noindent{\bf Remarks:}(1)~In the proof of Lemma \ref{RankLemma} we
could also have used \cite{Bi} for both the bundle case and the union
of two twisted $I$--bundles.  However, the setting of \cite{sou} is more
appropriate in this case (i.e. commensurable manifolds), and only a
mild extension is needed for us to handle the union of two twisted
$I$--bundles. Hence the reason for not using \cite{Bi} in this case. In
\S 4, we will need to use \cite{Bi}.\\[\baselineskip]

\noindent (2)~As noted in the introduction, we do not know if there
exists a pair of non-isometric, closed, orientable, hyperbolic
3--manifolds $M_1,M_2$ with $\mathcal{S}(M_1) =
\mathcal{S}(M_2)$. Since being either a virtual fiber or Fuchsian
depends only on $\Delta_\rho$ and not the ambient manifold, such a
pair would also satisfy both $\mathrm{S}_{vf}(M_1) =
\mathrm{S}_{vf}(M_2)$, $\mathrm{S}_{Fuc}(M_1) =
\mathrm{S}_{Fuc}(M_2)$. Examples where the latter equality holds were
constructed in \cite{MR} using a variation of the
method of Sunada \cite{Sun} for
constructing isospectral and length isospectral manifolds. That method
does not seem well-suited for also arranging equality between virtual
fibers. As with the full spectrum, we do not presently know if there
exists a pair of non-isometric, closed, hyperbolic 3--manifolds
$M_1,M_2$ with $\mathrm{S}_{vf}(M_1) = \mathrm{S}_{vf}(M_2)$.\\[\baselineskip]

\noindent (3)~All our results and methods of proof use the fact the hyperbolic 3-manifolds we consider are closed.  It would seem like an interesting problem to generalize the results of this paper to the case of finite volume non-compact hyperbolic 3-manifolds.

\section{A conjectural strengthening for arithmetic hyperbolic 3--manifolds}

In this section we deal with closed, arithmetic, hyperbolic 3--manifolds, and prove a stronger result (conjecturally) that involves only topological data. We refer the reader to \cite{MacR} for background on arithmetic hyperbolic 3--manifolds. Let us define {\em the topological virtual fiber set} of $M$ to be the set
\[ \mathrm{S}_{tvf}(M) = \{\Delta~:\Delta~\hbox{is isomorphic to a virtual fiber subgroup}\}. \]
Our strengthening relies on the following conjecture often referred to as the {\em short geodesic conjecture}.

\begin{conjecture}[Short Geodesic Conjecture]\label{short}
Let $M$ be a closed, orientable, arithmetic, hyperbolic 3--manifold.  Then there is a constant $C>0$ (independent of $M$) so that the length of the shortest closed geodesic in $M$ is at least $C$.
\end{conjecture}

Assuming this conjecture, we establish the following result.

\begin{theorem}\label{arithcase}
Assuming Conjecture \ref{short} there are at most finitely many closed orientable arithmetic hyperbolic 3--manifolds $M_1, M_2\ldots M_n$ so that $\mathrm{S}_{tvf}(M_i)=\mathrm{S}_{tvf}(M_j)$ for each $i,j$.
\end{theorem}

\noindent{\bf Proof:}~The proof of Theorem \ref{arithcase} is similar to the proof of Theorem \ref{mainprecise} and is done by contradiction.  If there is an infinite sequence of such manifolds $M_i$, by Borel \cite{Bo} their volumes are unbounded and Conjecture \ref{short} with the relationship between injectivity radius and systole implies that the injectivity radii are bounded from below. Choosing a minimal genus (topological) virtual fiber in each $M_i$ and applying Theorem \ref{soma}, it follows that for sufficiently large $i$, $M_i$ is either a genus $g$ fiber bundle or a union of two twisted $I$--bundles which is double covered by a genus $g$ fiber bundle. We now apply 
Biringer's extension of \cite{sou}, namely \cite[Thms 1.1, 5.2]{Bi}. That allows us to get control of the rank as in the proof of Lemma \ref{RankLemma}, and in particular, following the argument in the proof of Lemma \ref{RankLemma} leads to a similar contradiction on ranks as used in the proof of Theorem
\ref{mainprecise}.\qed


\noindent Department of Mathematics \\
Purdue University \\
West Lafayette, IN 47906 \\
email: {\tt dmcreyno@purdue.edu} \\

\noindent Department of Mathematics\\
Rice University\\
Houston, TX 77005\\
email: {\tt alan.reid@rice.edu}


\end{document}